\documentclass[12pt]{article}
\usepackage[english]{babel}
\usepackage{amsfonts}
\usepackage[dvips]{graphicx}

\begin{document}

\title{\bf A matrix generalization of Euler identity $e^{j\varphi}=cos\varphi + j\hspace{0.1cm}sin\varphi$}

\date{March 2007}

\author{\bf Gianluca Argentini \\
\normalsize gianluca.argentini@riellogroup.com \\
\textit{Research \& Development Department}\\
\textit{Riello Burners}, 37045 San Pietro di Legnago (Verona), Italy}

\maketitle

\begin{abstract}
In this work we present a matrix generalization of the Euler identity about exponential representation of a complex number. The concept of matrix exponential is used in a fundamental way. We define a notion of matrix imaginary unit which generalizes the usual complex imaginary unit. The Euler-like identity so obtained is compatible with the classical one. Also, we derive some exponential representation for matrix real and imaginary unit, and for the first Pauli matrix.
\end{abstract}

\section{The matrix exponential}

Let {\bf A} denotes a generic matrix. Based on the Taylor expansion centered at $0$ for the one real variable function $e^x$, the {\it matrix exponential} (see e.g. \cite{golub}) is formally defined as

\begin{equation}\label{matrixExp}
	e^{\bf A} = \sum_{n=0}^{+\infty}\frac{{\bf A}^n}{n!}
\end{equation}

\section{A class of complex matrices}

Let it be $j$ the imaginary unit, $j^2=-1$, $\alpha$, $a$, $b$ real numbers, $\alpha \neq 0$. We consider the class of 2$\times$2 matrices of the form

\[ {\bf T} = \left( \begin{array}{cc}\label{class}
a & jb \\
j\alpha^2b & a \end{array} \right)\]

\noindent These matrices are very important in the theory of transfer matrix method for modelization of acoustical transmission in physical structures (see \cite{munjal}, \cite{doria}). Note that $Det({\bf T})=a^2+\alpha^2b^2$, so that if $a$ and $b$ are not both zero, {\bf T} is invertible. In acoustical transfer matrices, usually $\alpha=\frac{S}{c}$, where $S$ is the section of a tube or duct and $c$ is the sound speed in the fluid contained in the tube.\\
Consider the matrix 

\[ {\bf \Phi} = \left( \begin{array}{cc}\label{imaginaryUnitMatrix}
0 & j \\
j\alpha^2 & 0 \end{array} \right)\]

\noindent Then, if ${\bf I}$ is the identity matrix, the following representation for previous ${\bf T}$ matrices holds:

\begin{equation}\label{representation}
	{\bf T} = a{\bf I} + b{\bf \Phi}
\end{equation}

\noindent Note the analogy with a usual complex number $a+jb$. The analogy is more evident if one consider that, with a simple calculation, ${\bf \Phi}^2 = -\alpha^2{\bf I}$. For this reason, we call ${\bf \Phi}$ the {\it imaginary unit matrix}. Also, note that, for $\alpha = 1$, we obtain ${\bf \Phi}=j\sigma_1$, where $\sigma_1$ is one of the {\it Pauli matrices} of quantum mechanics (see e.g. \cite{liboff}).

\section{A generalization of Euler identity}

\newtheorem{lemma}{Lemma}
\newtheorem{corol}{Corollary}
\newtheorem{theo}{Theorem}

The {\it Euler} identity $e^{j\varphi} = cos\varphi + j\hspace{0.1cm}sin\varphi$ is valid for any real $\varphi$. Usually this formula is proven by use of Taylor expansion of the complex function $e^z$ and of the real functions $cos\varphi$ and $sin\varphi$ (see \cite{boyer}).\\
We prove a generalization, in the environment of the complex matrices of type (\ref{representation}), of this identity.

\begin{lemma}
	Let it be ${\bf \Psi} = -j{\bf \Phi}$. Then, for every natural $n$, the following relation holds:\\
	\begin{equation}\label{psiProp}
		{\bf \Psi}^n = \alpha^{n-r}{\bf \Psi}^r
	\end{equation}
where $r = mod(n,2)$.
\end{lemma}
	
\noindent {\it Dim}. By induction on $n$. For $n=0$ and $n=1$ the thesis is obvious. Note that

\[ {\bf \Psi} = \left( \begin{array}{cc}
0 & 1 \\
\alpha^2 & 0 \end{array} \right)\]

\noindent Let it be $n=2$. By a simple calculation

\[ {\bf \Psi}^2 = \left( \begin{array}{cc}
0 & 1 \\
\alpha^2 & 0 \end{array} \right)^2 = \left( \begin{array}{cc}
\alpha^2 & 0 \\
0 & \alpha^2 \end{array} \right) = \alpha^2{\bf I} = \alpha^2{\bf \Psi}^0 \]

\noindent and the thesis is verified. Then we suppose that the thesis is verified for a generic $n$ too. Using the last inductive step we have

\begin{equation}
	{\bf \Psi}^{n+1} = {\bf \Psi}^n{\bf \Psi} = \alpha^{n-r}{\bf \Psi}^r{\bf \Psi}
\end{equation}

\noindent where $r = mod(n,2)$. If $mod(n+1,2)=1$, then $r=0$, therefore

\begin{equation}
	{\bf \Psi}^{n+1} = \alpha^n{\bf \Psi} = \alpha^{n+1-s}{\bf \Psi}^s
\end{equation}

\noindent with $s = mod(n+1,2) = 1$, and the thesis is true in this case. If $mod(n+1,2)=0$, then $r=1$, therefore, using the first inductive step for $n=2$,
	
\begin{equation}
	{\bf \Psi}^{n+1} = \alpha^{n-1}{\bf \Psi}{\bf \Psi} = \alpha^{n+1}{\bf \Psi}^0 = \alpha^{n+1-s}{\bf \Psi}^s
\end{equation}

\noindent with $s = mod(n+1,2) = 0$, and the thesis is true in this case too. $\square$\\

Now we can prove the matrix generalization of Euler identity:

\begin{theo}
	For every real $\varphi$, the following formula holds:
	\begin{equation}\label{genEuler}
		e^{\varphi {\bf \Phi}} = cos(\alpha \varphi){\bf I} + \frac{1}{\alpha}sin(\alpha \varphi){\bf \Phi}
	\end{equation}
\end{theo}

\noindent {\it Dim}. If $\varphi = 0$ the formula is obvious. For $\varphi \neq 0$, from the formal definition (\ref{matrixExp}) we can write

\begin{equation}
	e^{\varphi {\bf \Phi}} = \sum_{n=0}^{+\infty}\frac{(\varphi {\bf \Phi})^n}{n!} = \sum_{n \hspace{0.1cm} odd}^{+\infty}\frac{\varphi^n {\bf \Phi}^n}{n!}+\sum_{n \hspace{0.1cm} even}^{+\infty}\frac{\varphi^n {\bf \Phi}^n}{n!}
\end{equation}

\noindent Recall that, if $n$ is even, then $j^n$ alternates $-1$ and $+1$, while if $n$ is odd, then $j^n$ alternates $-j$ and $+j$. Therefore, from ${\bf \Phi}=j{\bf \Psi}$, from the usual series expansion for $cos(\alpha \varphi)$ and $sin(\alpha \varphi)$, and using the previous Lemma, we have

\begin{eqnarray}
	e^{\varphi {\bf \Phi}} = \sum_{n \hspace{0.1cm} even}^{+\infty}j^n\frac{\varphi^n {\bf \Psi}^n}{n!} + \sum_{n \hspace{0.1cm} odd}^{+\infty}j^n\frac{\varphi^n {\bf \Psi}^n}{n!} =\\ \nonumber
	= \left(\sum_{n \hspace{0.1cm} even}^{+\infty}j^n\frac{(\alpha \varphi)^n}{n!}\right){\bf I} + \frac{1}{\alpha} \left(\sum_{n \hspace{0.1cm} odd}^{+\infty}j^n\frac{(\alpha \varphi)^n}{n!}\right){\bf \Psi} = \\ \nonumber
	= cos(\alpha \varphi){\bf I} + \frac{1}{\alpha}j\hspace{0.1cm}sin(\alpha \varphi){\bf \Psi} = cos(\alpha \varphi){\bf I} + \frac{1}{\alpha}sin(\alpha \varphi){\bf \Phi} 
\end{eqnarray}

\noindent that is the thesis. $\square$\\

\noindent {\bf Note 1}. Let it be $\alpha=1$, and ${\bf I}=[1]$, ${\bf \Phi}=[j]$ two $1\times1$ {\it matrices}, so that ${\bf I}$ is the usual real unit and ${\bf \Phi}$ the usual imaginary unit. From (\ref{genEuler}) we have

\begin{equation}
	e^{j\varphi} = cos(\varphi)[1]+ sin(\varphi)[j] = cos\varphi+j\hspace{0.1cm}sin\varphi
\end{equation}

\noindent that is the classical Euler identity.\\

\noindent {\bf Note 2}. If we write in explicite mode the relation (\ref{genEuler}), we obtain

\[ e^{\varphi {\bf \Phi}} = \left( \begin{array}{cc}
cos(\alpha \varphi) & j\frac{1}{\alpha}sin(\alpha \varphi) \\
\null & \null \\
j\alpha\hspace{0.1cm}sin(\alpha \varphi) & cos(\alpha \varphi) \end{array} \right)\]

\noindent so that $Det(e^{\varphi {\bf \Phi}}) = 1$, which is compatible with the fact that for usual complex numbers $\left|e^{j\varphi}\right|=1$.\\

\noindent {\bf Note 3}. If ${\alpha=1}$ and ${\varphi=2m\pi}$, with $m$ integer, (\ref{genEuler}) becomes

\begin{equation}
	e^{2m\pi {\bf \Phi}} = {\bf I}
\end{equation}

\noindent that is a matrix unit representation. The classical analogous formula is $e^{j2m\pi}=1$.\\

\noindent {\bf Note 4}. If ${\alpha=1}$ and ${\varphi=m\frac{\pi}{2}}$, with $m=1+4k$, $k$ integer, (\ref{genEuler}) becomes

\begin{equation}
	e^{m\frac{\pi}{2} {\bf \Phi}} = {\bf \Phi}
\end{equation}

\noindent that is a matrix imaginary unit representation. The classical analogous formula is $e^{jm\frac{\pi}{2}}=j$. Also, if we multiply previous formula by $-j$, we have an exponential representation of Pauli matrix ${\sigma_1}$:

\begin{equation}
	{\sigma_1} = -je^{m\frac{\pi}{2} {\bf \Phi}}
\end{equation}

\end{document}